\documentclass{article}
\usepackage{amsmath}
\usepackage{amsfonts}

\newtheorem{theorem}{Theorem}

\newtheorem{definition}[theorem]{Definition}
\newtheorem{example}[theorem]{Example}

\newtheorem{lemma}[theorem]{Lemma}

\newenvironment{proof}[1][Proof]{\noindent\textbf{#1.} }{\ \rule{0.5em}{0.5em}}
\input{tcilatex}
\begin{document}

\title{Pattern Avoiding Ballot Paths and Finite Operator Calculus}
\author{Heinrich Niederhausen and Shaun Sullivan \\
Florida Atlantic University, Boca Raton, Florida}
\maketitle

\begin{abstract}
Counting pattern avoiding ballot paths begins with a careful analysis of the
pattern. Not the length, but the characteristics of the pattern are
responsible for the difficulties in finding explicit solutions. Certain
features, like overlap and difference in number of $\rightarrow $ and $%
\uparrow $ steps determine the solution of the recursion formula. If the
recursion can be solved by a polynomial sequence, we apply the Finite
Operator Calculus to find an explicit form of the solution in terms of
binomial coefficients.

\end{abstract}

\section{Introduction}

A ballot path stays weakly above the diagonal $y=x$, starts at the origin
and takes steps from $\{\uparrow ,\rightarrow \}$. A pattern is a finite
string made from the same step set; it is also a path. Dyck paths are
equivalent to ballot paths, taking steps from $\left\{ \nearrow ,\searrow
\right\} $, staying weakly above the $x$-axis. Dyck paths containing $k$
strings of length $3$ were discussed by E. Deutsch in \cite{Deutsch}. One of
the most recent papers on patterns of length $4$ occurring $k$ times in Dyck
paths was written by A. Sapounakis, I. Tasoulas, P. Tsikouras, Counting
strings in Dyck paths, 2007, to appear in \emph{Discrete Mathematics} \cite%
{stt}.\textit{\ }The authors find generating functions for all $16$ patterns
in Dyck path returning to the $x$-axis. Returning to the $x$-axis at the end
of the path has the advantage that going backwards through a path we find
the reversed pattern exactly the same number of times (see Section \ref%
{SecRightSteps}). This reduces significantly the number of patterns under
consideration.

For easier presentation we decided to talk about ballot instead of Dyck
paths. We will look at paths that end at any point $\left( n,m\right) $
above or on the main diagonal. That will give us more cases to consider,
because the path reversal bijection will be of no help in general. On the
other hand, we will investigate \emph{pattern avoiding paths} only ($k=0$),
but we will do this for patterns of general length. We use the Finite
Operator Calculus approach, as developed by G.-C. Rota, D. Kahaner, and A.
Odlyzko \cite{ROK}. An interpretation of this theory in view of solving
recursions is given in \cite{Nied03}. The Finite Operator Calculus works
with polynomials; it applies to pattern avoiding ballot paths, because for
many patterns we can find a polynomial sequence whose values enumerate all
possible paths.

For convenience, we let $u=\uparrow $ and $r=\rightarrow $. Of course,
avoiding the pattern $p$ may imply avoiding a pattern that contains $p$
twice in some overlapping form, like $uruuru$ is contained twice in $uruur%
\mathbf{u}ruuru$, but also twice in $uru\mathbf{uru}uru$. This problem of
overlaps makes it harder to find recursions. There are two properties of the
patterns we avoid that determine the recurrence relations.

\begin{definition}
The bifix index of a pattern $p$ is the number of distinct nonempty patterns 
$o$ such that $p$ that can be written in the form $p=op^{\prime }$ and $%
p=p^{\prime \prime }o$ for $o,p^{\prime },p^{\prime \prime }\in
\{u,r\}^{\ast }$. If a pattern has bifix index 0, then we call it bifix-free.
\end{definition}

The example above, $uruuru$, has bifix index $2$. We restrict ourselves to
bifix index $0$ or $1$ in this paper. However, there is a notable exception;
a pattern of only right steps has a high bifix index. We can avoid this
pattern with a prefixed up step, which makes the prefixed pattern bifix
free. Another exception is a pattern of only up steps. We will begin in
Section \ref{SecRightSteps} with this case to explain our approach.

\begin{definition}
\label{DefDepth}Let $d(p)$ be the number of $u$'s minus the number of $r$'s
in the pattern $p$. The depth of $p$ is $\max \{d(p^{\prime })\mid
p=qp^{\prime },q\in \{u,r\}^{\ast }\}$.
\end{definition}

Intuitively, the depth is the maximum distance $p$ goes below the line $y=x$
when the end of $p$ is attached to the line $y=x$. Notice that since the
empy pattern $\epsilon $ is in $\{p^{\prime }|p=qp^{\prime }\}$, the depth
of a pattern is always a nonnegative integer.

If $a$ is the number of $r$'s in $p$ and $c$ is the number of $u$'s, then we
say $p$ has dimensions $a\times c$.

We consider patterns of the following form:

\begin{enumerate}
\item $r^{a}$ and $u^{c}$ for any $a,c\geq 0$.

\item Bifix-free patterns of depth $0$ and dimensions $a\times c$, with $%
a\geq c\geq 1$, and $a\geq 2$.

\item Patterns of depth $0$ with bifix index 0 or 1 (with the same
restrictions as above, and corresponding restriction about the $op^{\prime }$
piece of $p=op^{\prime }o$).

\item Patterns with depth at least $1$ and bifix index $0$ or $1$ (with
similar restrictions).

\item All patterns of length 4.
\end{enumerate}

\section{Lattice paths avoiding $r^{a}$ and $u^{c}$\label{SecRightSteps}}

\begin{definition}
Let $s_{n}(m;p,l)$ be the number of lattice paths from $(0,0)$ to $(n,m)$
staying weakly above the line $y=x-l$ avoiding the pattern $p\in \{\uparrow
,\rightarrow \}^{\ast }$.
\end{definition}

Note that if $l=0$ we count pattern avoiding\emph{\ }ballot paths. A ballot
path must start with an up step; if the pattern to avoid is $r^{a}$, we are
really avoiding $ur^{a}$. In other words, $s_{n}\left( m;r^{a},0\right)
=s_{n}\left( m;ur^{a},0\right) $. This means that the bifix-rich pattern $%
r^{a}$ can be replaced by the bifix-free pattern $ur^{a}$. We loose this
property if $l>0$. However, the recursion%
\begin{equation*}
s_{n}(m;r^{a},l)=s_{n-1}(m;r^{a},l)+s_{n}(m-1;r^{a},l)-s_{n-a}(m-1;r^{a},l)
\end{equation*}%
still holds for all $m\geq n>0$. Only the initial values change with $l$, $%
s_{n}\left( 0;r^{a},l\right) =1$ for $0\leq n\leq \min \left( a-1,l\right) $%
, and $0$ else. The recursion says that we must subtract from the ballot
recursion the paths reaching $(n,m)$ with \textit{exactly} $a$ down steps at
the end, and one up step in the beginning. These are counted by $%
s_{n-a}(m-1;r^{a},l)$.

\begin{equation*}
\begin{tabular}{l}
$%
\begin{tabular}{r||rrrrrrrrrr}
$m$ &  &  &  &  &  &  &  &  &  &  \\ 
$4$ & 1 & 5 & 15 & 35 & 65 & 101 & 135 & 155 & 152 & 112 \\ 
$3$ & 1 & 4 & 10 & 20 & 31 & 40 & 44 & 40 & 28 & 0 \\ 
$2$ & 1 & 3 & 6 & 10 & 12 & 12 & 10 & 6 & 0 & -16 \\ 
$1$ & 1 & 2 & 3 & 4 & 3 & 2 & 1 & 0 & -3 & -14 \\ 
$0$ & 1 & 1 & 1 & 1 & 0 & 0 & 0 & 0 & -3 & -11 \\ \hline\hline
$n:$ & $0$ & $1$ & $2$ & $3$ & $4$ & $5$ & $6$ & $7$ & $8$ & $9$%
\end{tabular}%
$ \\ 
\multicolumn{1}{c}{$s_{n}\left( m;r^{4},5\right) $}%
\end{tabular}%
\end{equation*}%
The table shows the number of paths avoiding $r^{4}$ or $rrrr$ above the
zeros and extending below the boundary using the recurrence to obtain a
polynomial sequence. 

Induction over $n$ shows that $(s_{n})$ is a polynomial sequence with $\deg {%
s_{n}}=n$. Using operators on polynomials, we can write the recurrence as

\begin{equation*}
1-E^{-1}=B-B^{a}E^{-1}
\end{equation*}%
where the linear operators $B$ and $E^{v}$ are defined by linear extension
of $Bs_{n}(x)=s_{n-1}(x)$ and $E^{v}s_{n}(x)=s_{n}(x+v)$, the shift by $v$.
The operators $\nabla =1-E^{-1}$ and $E^{-1}$ both have power series
expansions in $D$, the derivative operator. Such operators are called \emph{%
shift-invariant}. Hence $B$ must be shift-invariant too, and therefore
commutes with $\nabla $ and $E^{v}$. The power series for $B$ must be of
order $1$, because $B$ reduces degrees by $1$. Such linear operators are
called \textit{delta operators}. The basic sequence $\left( b_{n}\left(
x\right) \right) _{n\geq 0}$ of a delta operator $B$ is a sequence of
polynomials such that $\deg b_{n}=n$, $Bb_{n}\left( x\right) =b_{n-1}\left(
x\right) $ (like the \textit{Sheffer sequence} $s_{n}\left( x\right) $ for $B
$), and initial conditions $b_{n}\left( 0\right) =\delta _{0,n}$ for all $%
n\in \mathbb{N}_{0}$. In our special case, the basic sequence is easily
determined. Solving for $E^{1}$ shows that

\begin{equation*}
E^{1}=\sum\limits_{i=0}^{a-1}B^{i}.
\end{equation*}

Finite Operator Calculus tells us \cite[(2.5)]{Nied03} that if $%
E^{1}=1+\sigma (B)$, where $\sigma (t)$ is a power series of order 1, then
the basic sequence $b_{n}(x)$ of $B$ has the generating function

\begin{equation*}
\sum\limits_{n\geq 0}b_{n}(x)t^{n}=(1+\sigma (t))^{x}.
\end{equation*}%
Thus, in our case $b_{n}(x)=\left[ t^{n}\right] (1+t+t^{2}+\cdots
+t^{a-1})^{x}$

\begin{definition}
The \emph{geometric coefficient }is defined as

\begin{equation*}
\dbinom{x}{n}_{a}=[t^{n}](1+t+\cdots
+t^{a-1})^{x}=\sum\limits_{i=0}^{\left\lfloor n/a\right\rfloor }(-1)^{i}%
\dbinom{x}{i}\dbinom{x+n-ai-1}{n-ai}
\end{equation*}
\end{definition}

Note that for $a=2$ the geometric coefficient equals the binomial
coefficient $\dbinom{x}{n}$. These numbers have already been studied by
Euler \cite{Euler}. Therefore they are also called \emph{Eulerian
coefficients}. Some interesting properties of geomeric coefficients are
given in \cite{N-Su07}.

Thus $b_{n}\left( x\right) =\dbinom{x}{n}_{a}$. For\emph{\ }$l=0$ the
Sheffer sequence $(s_{n})$ has initial values $s_{n}(n-1)=\delta _{n,0}$.
Abelization \cite{Nied03} gives us

\begin{equation*}
s_{n}(x;r^{a},0)=\dfrac{x-n+1}{x+1}b_{n}(x+1)=\dfrac{x-n+1}{x+1}\dbinom{x+1}{%
n}_{a}
\end{equation*}%
so the number of ballot paths avoiding $r^{a}$ is

\begin{equation}
s_{n}(n;r^{a},0)=\dfrac{1}{n+1}\dbinom{n+1}{n}_{a}.  \label{(l=0)}
\end{equation}%
The reflection of a pattern $p$\ is the pattern where every up step becomes
a right step, and vice versa. The pattern $\tilde{p}$ is the \emph{reverse}
of the pattern $p$, if it is the reflection of $p$ read backwards; for
example if $p=uruuruu$, then $\tilde{p}=rrurrur$. It is a fundamental
principal in pattern avoidance, proved by reflection, that the number of $%
\{\uparrow ,\rightarrow \}$ lattice paths weakly above $y=x$ ending at $%
\left( n,n+l\right) $ avoiding $p$ equals the number of $\{\uparrow
,\rightarrow \}$ lattice paths weakly above $y=x-l$ ending at $\left(
n+l,l\right) $ avoiding the reverse pattern $\tilde{p}$,%
\begin{equation*}
s_{n}(n+l;p,0)=s_{n+l}(n;\tilde{p},l).
\end{equation*}%
If we are only interested in paths returning to the diagonal ($l=0$), we can
see how this principal saves us a great deal of work. For general $l>0$,
this will not be the case, because we will not be able to find the paths
avoiding $\tilde{p}$ staying weakly above $y=x-l$. There is one notable
exception, $p=r^{a}$. For general $l$, the initial values $s_{n}\left(
0;r^{a},l\right) $ agree up to $n=l$ with $b_{n}\left( 1\right) $, and
therefore $s_{n}\left( x;r^{a},l\right) =b_{n}\left( x+1\right) $ for all $%
n\leq l$. For $n>l$ we have $s_{n}\left( n-l-1;r^{a},l\right) =0$. The \emph{%
Binomial Theorem for Sheffer Sequences}\ \cite{Nied03} expands $\left(
s_{n}\right) $ under these and similar initial values.

\begin{theorem}
If $(t_{n})$ is a Sheffer sequence and $(q_{n})$ the basic sequence for the
same delta operator, then%
\begin{equation*}
t_{n}(x+y)=\sum\limits_{k=0}^{n}t_{k}(y+vk)\dfrac{x-vn}{x-vk}q_{n-k}(x-vk)
\end{equation*}%
for all $v\in \mathbb{R}$.
\end{theorem}

Hence%
\begin{eqnarray}
s_{n+l}(n;r^{a},l) &=&\sum\limits_{k=0}^{n+l}s_{k}(k-l-1;r^{a},l)\dfrac{1}{%
n+l+1-k}b_{n+l-k}(n+l+1-k)  \notag \\
&=&\sum\limits_{k=0}^{l}\dfrac{1}{n+l+1-k}\dbinom{k-l}{k}_{a}\dbinom{n+l+1-k%
}{n+l-k}_{a}.  \notag
\end{eqnarray}%
Of course, this expansion reduces to (\ref{(l=0)}) if $l=0$.

Because we were able to find the number of paths weakly above $y=x-l$ for
general $l$ as polynomials, we are also able to apply the above general
principle, saying that $s_{n}(n+l;u^{c},0)=s_{n+l}\left( n;r^{c},l\right) $.
Thus%
\begin{equation}
s_{n}(m;u^{c},0)=\sum\limits_{k=0}^{m-n}\dfrac{1}{m+1-k}\dbinom{k-m+n}{k}_{c}%
\dbinom{m+1-k}{m-k}_{c}  \notag
\end{equation}%
\begin{equation*}
\begin{tabular}{c}
$%
\begin{tabular}{r||rrrrr}
$9$ &  &  & 1 & 19 & \textbf{112} \\ 
$8$ &  &  & 3 & \textbf{28} & 116 \\ 
$7$ &  &  & \textbf{6} & 33 & 101 \\ 
$6$ &  & \textbf{1} & 9 & 32 & 68 \\ 
$5$ & \textbf{0} & 2 & 10 & 23 & 36 \\ 
$4$ &  & 3 & 8 & 13 & 13 \\ 
$3$ & 1 & 3 & 5 & 5 &  \\ 
$2$ & 1 & 2 & 2 &  &  \\ 
$1$ & 1 & 1 &  &  &  \\ 
$0$ & 1 &  &  &  &  \\ \hline\hline
$n:$ & $0$ & $1$ & $2$ & $3$ & $4$%
\end{tabular}%
$ \\ 
$s_{n}\left( m;u^{4},0\right) $ (the numbers in bold agree with $%
s_{n+5}\left( n;r^{4},5\right) $ in the Table above)%
\end{tabular}%
\end{equation*}

\section{Bifix-free patterns with depth 0}

From now on we will only look at paths weakly above the diagonal $y=x$. We
write $s_{n}\left( x;p\right) $\ instead of $s_{n}\left( x;p,0\right) $, and
we also may omit the pattern $p$\ from the notation. As with the pattern $%
r^{a}$, we can find a single recurrence relation that holds everywhere in
the octant. If the pattern is bifix-free, we need only to subtract paths
that would end in the pattern, thus we have the recurrence

\begin{equation}
s_{n}(m;p)=s_{n-1}(m;p)+s_{n}(m-1;p)-s_{n-a}(m-c;p)  \label{(bifix0)}
\end{equation}%
where $p$ has dimensions $a\times c$. For example $uurrurrur$ has dimensions 
$5\times 4$, and depth $0$. The recurrence has a polynomial solution if the
depth is $0$,\ and $a\geq c\geq 1$, $a\geq 2$ (If $a=1$ then $p=ur$, a
pattern we do not want to avoid). In operators:

\begin{equation*}
\nabla=B(1-B^{a-1}E^{-c})
\end{equation*}

Since the delta operator $\nabla $ can be written as a delta series in $B$,
the operator $B$ is also a delta operator. The basic sequence can be
expressed via the Transfer Formula \cite[Theorem 1]{Nied03}:

\begin{equation*}
b_n(x)=x\sum\limits_{i=0}^{\frac{n}{a-1}}\frac{(-1)^i}{x-ci}\dbinom{n-(a-1)i%
}{i}\dbinom{x+n-(a+c-1)i-1}{n-(a-1)i}
\end{equation*}

Since our initial values are $s_{n}(n-1;p)=\delta _{n,0}$, we use
Abelization \cite[(2.5)]{Nied03} to obtain:

\begin{equation*}
s_n(x)=(x-n+1)\sum\limits_{i=0}^{\frac{n}{a-1}}\frac{(-1)^i}{x-ci+1}\dbinom{%
n-(a-1)i}{i}\dbinom{x+n-(a+c-1)i}{n-(a-1)i}
\end{equation*}

Therefore the number of\ ballot paths avoiding $p$ and returning to the
diagonal is

\begin{equation*}
s_n(n)=\sum\limits_{i=0}^{\frac{n}{a-1}}\frac{(-1)^i}{n-ci+1}\dbinom{n-(a-1)i%
}{i}\dbinom{2n-(a+c-1)i}{n-(a-1)i}
\end{equation*}

\section{Patterns with depth 0 and bifix index 1}

The pattern $p$ has bifix index $1$ if there exists a unique nonempty
pattern $o$ such that $p=op^{\prime }o$. If $p$ has dimensions $a\times c$
and $op^{\prime }$ has dimensions $b\times d$, we have a recurrence of the
form

\begin{equation}
s_{n}(x)=s_{n-1}(x)+s_{n}(x-1)-\sum\limits_{i\geq
0}(-1)^{i}s_{n-a-bi}(x-c-di)  \label{(bifix1)}
\end{equation}

For example let $p=urruurr$. It has depth $0$, with dimensions $4\times 3$
and bifix $urr$, so $b=2$ and $d=2$. From the paths reaching $(n-a,x-c)$
those ending in $urru$ cannot be included in the recurrence and must be
subtracted, and from those again we cannot include paths ending in $urru$,
and so on. The $op^{\prime }$\ piece of the pattern that is responsible for
this exclusion-inclusion process may not go below the diagonal; hence it
must have depth\emph{\ }$\leq a-c$\emph{.}

From examining the summation in the recurrence, we notice that we must have $%
b\geq d$. If $b<d$, then at some point we would be using numbers below the $%
y=x$ boundary, which are only the polynomial extensions and do not count
paths. Notice in our above example $b-d=0$, thus the recurrence formula
applies to the pattern $urruurr$. 
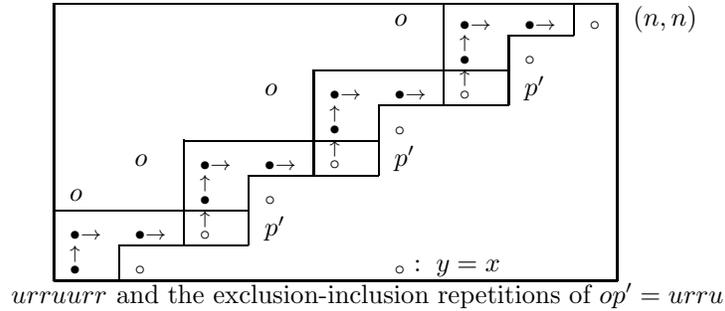
\begin{equation*}
\begin{tabular}{c}
$%
\begin{tabular}{l|lllllllll|l}
\cline{2-10}\cline{9-9}
&  &  &  &  &  & $o$ & \multicolumn{1}{|l}{$_{\bullet \rightarrow }$} & $%
_{\bullet \rightarrow }$ & \multicolumn{1}{|l|}{$_{\circ }$} & $\left(
n,n\right) $ \\ \cline{9-9}
&  &  &  &  &  &  & \multicolumn{1}{|l|}{$_{\bullet }^{\uparrow }$} & 
\multicolumn{1}{|l}{$_{\circ }$} &  &  \\ \cline{6-8}
&  &  &  & $o$ & \multicolumn{1}{|l}{$_{\bullet \rightarrow }$} & $_{\bullet
\rightarrow }$ & \multicolumn{1}{|l}{$_{\circ }^{\uparrow }$} & 
\multicolumn{1}{|l}{$p^{\prime }$} &  &  \\ \cline{7-8}
&  &  &  &  & \multicolumn{1}{|l}{$_{\bullet }^{\uparrow }$} & 
\multicolumn{1}{|l}{$_{\circ }$} &  &  &  &  \\ \cline{4-6}
&  & $o$ & \multicolumn{1}{|l}{$_{\bullet \rightarrow }$} & $_{\bullet
\rightarrow }$ & \multicolumn{1}{|l}{$_{\circ }^{\uparrow }$} & 
\multicolumn{1}{|l}{$p^{\prime }$} &  &  &  &  \\ \cline{5-6}
& $o$ &  & \multicolumn{1}{|l}{$_{\bullet }^{\uparrow }$} & 
\multicolumn{1}{|l}{$_{\circ }$} &  &  &  &  &  &  \\ \cline{2-4}
& $_{\bullet \rightarrow }$ & $_{\bullet \rightarrow }$ & 
\multicolumn{1}{|l}{$_{\circ }^{\uparrow }$} & \multicolumn{1}{|l}{$%
p^{\prime }$} &  & \multicolumn{4}{l|}{} &  \\ \cline{3-4}
& $_{\bullet }^{\uparrow }$ & \multicolumn{1}{|l}{$_{\circ }$} &  &  &  & 
\multicolumn{4}{l|}{$_{\circ }:\;y=x$} &  \\ \cline{2-10}
\end{tabular}%
$ \\ 
$urruurr$ and the exclusion-inclusion repetitions of $op^{\prime }=urru$%
\end{tabular}%
\end{equation*}

Summary of the conditions on $p$:

\begin{enumerate}
\item $p=op^{\prime }o$, where $o$ is nonempty and unique (bifix index $=1$).

\item $p$ has dimensions $a\times c$, where $a\geq c\geq 1$, $a\geq 1$.

\item $p^{\prime }o$ has dimensions $b\times d$, where $b\geq d$, $b\geq 1$.

\item $\func{depth}\left( p^{\prime }o\right) \leq a-c$.
\end{enumerate}

We assume that $b\geq 1$, because otherwise $p=u^{c}$, a case we considered
already in a previous section. Under these conditions $\left( s_{n}\right) $
is a Sheffer sequence, and the recursion (\ref{(bifix1)}) can be written in
terms of the operator $B:s_{n}\mapsto s_{n-1}$ as%
\begin{equation*}
1=B+E^{-1}-\sum\limits_{i\geq 0}(-1)^{i}B^{a+bi}E^{-c-di}=B+E^{-1}-\dfrac{%
B^{a}E^{-c}}{1+B^{b}E^{-d}}
\end{equation*}%
or

\begin{eqnarray}
\nabla &=&B+B^{b+1}E^{-d}+B^{b}E^{-d-1}-B^{a}E^{-c}-B^{b}E^{-d}
\label{(index1)} \\
&=&B-B^{a}E^{-c}+B^{b+1}E^{-d}-B^{b}E^{-d}\nabla  \notag
\end{eqnarray}

We view $\nabla $ as a formal power series in $B$ with coefficients in the
ring of shift-invariant operators. The generalized Transfer Formula \cite[%
Theorem 2]{Nied03} tells us that%
\begin{equation*}
b_{n}\left( x\right) =x\sum_{i=1}^{n}C_{n,i}\frac{1}{i}\binom{i-1+x}{i-1}
\end{equation*}%
is the basic polynomial for $B$, where 
\begin{equation*}
C_{n,i}=\left[ B^{n}\right] \left(
B-B^{a}E^{-c}+B^{b+1}E^{-d}-B^{b}E^{-d}\nabla \right) ^{i}
\end{equation*}%
Hence%
\begin{eqnarray*}
b_{n}\left( x\right)  &=&x\sum_{j,k,l\geq 0:\ \left( a-1\right) j+bk+\left(
b-1\right) l\leq n-1}\binom{n-\left( a-1\right) j-bk-\left( b-1\right) l}{%
j,k,l}\left( -1\right) ^{j+l} \\
&&\times E^{-cj-dk-dl}\nabla ^{l}\frac{1}{n-\left( a-1\right) j-bk-\left(
b-1\right) l}\binom{n-\left( a-1\right) j-bk-\left( b-1\right) l-1+x}{%
n-\left( a-1\right) j-bk-\left( b-1\right) l-1} \\
&=&x\sum_{j=0}^{\frac{n-1}{a-1}}\sum_{k=0}^{\frac{n-1-\left( a-1\right) j}{b}%
}\sum_{l=0}^{\frac{n-1-\left( a-1\right) j}{b}-k}\binom{n-\left( a-1\right)
j-bk-\left( b-1\right) l}{j,k,l} \\
&&\times \frac{\left( -1\right) ^{j+l}}{n-\left( a-1\right) j-bk-\left(
b-1\right) l}\binom{n-\left( a+c-1\right) j-\left( b+d\right) \left(
k+l\right) -1+x}{n-\left( a-1\right) j-b\left( k+l\right) -1}
\end{eqnarray*}%
Because $s_{n}\left( n-1\right) =\delta _{0,n}$ we still get $s_{n}\left(
x\right) =\frac{x-n+1}{x+1}b_{n}\left( x+1\right) $.

\begin{equation*}
\begin{tabular}{l}
$%
\begin{tabular}{l||llllllllll}
$m$ & 1 & 8 & 35 & 110 & 270 & 544 & 920 & 1272 & 1236 & 0 \\ 
$7$ & 1 & 7 & 27 & 75 & 161 & 279 & 389 & 377 & 0 &  \\ 
$6$ & 1 & 6 & 20 & 48 & 87 & 122 & 118 & 0 &  &  \\ 
$5$ & 1 & 5 & 14 & 28 & 40 & 38 & 0 &  &  &  \\ 
$4$ & 1 & 4 & 9 & 14 & 13 & 0 &  &  &  &  \\ 
$3$ & 1 & 3 & 5 & 5 & 0 &  &  &  &  &  \\ 
$2$ & 1 & 2 & 2 & 0 &  &  &  &  &  &  \\ 
$1$ & 1 & 1 & 0 &  &  &  &  &  &  &  \\ 
$0$ & 1 & 0 &  &  &  &  &  &  &  &  \\ \hline\hline
& $0$ & $1$ & $2$ & $3$ & $4$ & $5$ & $6$ & $7$ & $8$ & $n$%
\end{tabular}%
$ \\ 
The number of ballot paths avoiding $urruurr$%
\end{tabular}%
\end{equation*}

\begin{example}
For the pattern $urruurr$ we find the parameters $a=\#r$'s$\ =4$, $c=\#u$'s$%
~=3$, $b=2$, $d=2$. Because the $\func{depth}\left( urruurr\right) =0$, we
obtain the explicit solution\newline
$s\left( n,x\right) =\sum_{j=0}^{\left\lfloor \left( n-1\right)
/3\right\rfloor }\sum_{k=0}^{\left\lfloor \left( n-1-3j\right)
/2\right\rfloor }\sum_{l=0}^{\left\lfloor \left( n-1-3j-2k\right)
/2\right\rfloor }\binom{n-3j-2k-l}{j}\binom{n-4j-2k-l}{k}\binom{n-4j-3k-l}{l}%
\times $\newline
$\qquad \qquad \times \frac{\left( x-n+1\right) \left( -1\right) ^{j+l}}{%
n-3j-2k-l}\binom{n-6j-4\left( k+l\right) +x}{n-3j-2\left( k+l\right) -1}$
\end{example}

The above operator equation (\ref{(index1)}) simplifies when $a=b+1$ and $%
c=d $. This corresponds to a pattern of the form $rp^{\prime }r$.\emph{\ }%
For this case we get%
\begin{equation*}
\nabla =B(1+B^{b}E^{-d})^{-1}
\end{equation*}%
\newline
and%
\begin{equation*}
s_{n}(x)=\frac{(x-n+1}{x+1}b_{n}\left( x+1\right)
=\sum\limits_{i=0}^{\left\lfloor n/b\right\rfloor }\dfrac{(-1)^{i}}{x-di+1}%
\dbinom{n-(b-1)i-1}{i}\dbinom{x+n-(d+b+1)i}{n-bi}
\end{equation*}%
The number of paths returning to the diagonal equals

\begin{equation*}
s_{n}(n)=\sum\limits_{i=0}^{\left\lfloor n/b\right\rfloor }\dfrac{(-1)^{i}}{%
n-di+1}\dbinom{n-(b-1)i-1}{i}\dbinom{2n-(d+b+1)i}{n-bi}.
\end{equation*}

\section{Patterns with depth at least 1 and bifix index 0 or 1}

Up to now we studied patterns such that the pattern avoiding ballot numbers
could be continued to Sheffer polynomials below the diagonal, using the
notation $s_{n}\left( x\right) $. Because that is no longer true in this
Section, we begin with the new notation $D\left( n,m\right) $\ for the
pattern avoiding ballot numbers reaching $\left( n,m\right) $. It is easy to
see that if the pattern has depth $\delta >0$, then the recurrence for all
points $\left( n,m\right) $\emph{\ }between the lines $y=x$ and $y=x+\delta $
is

\begin{equation*}
D(n,m)=D(n-1,m)+D(n,m-1)
\end{equation*}%
since by definition of depth a path ending at $(n,m)$ with $n+\delta >m\geq n
$ cannot contain the pattern $p$. In general, these numbers can not be
extended to polynomials of degree $n$. As initial values we have $D\left(
n,n-1\right) =\delta _{0,n}$ for all $n\geq 0$.

If $\left( n,m\right) $ falls weakly\emph{\ }above the line $y=x+\delta $ we
have the recurrence%
\begin{equation}
D(n,m)=D(n-1,m)+D(n,m-1)-D(n-a,m-c)  \label{(bifix0')}
\end{equation}%
for bifix index 0, and

\begin{equation}
D(n,m)=D(n-1,m)+D(n,m-1)-\sum\limits_{i\geq 0}(-1)^{i}D(n-a-bi,m-c-di)
\label{(bifix1')}
\end{equation}%
for bifix index 1 (see (\ref{(bifix0)}) and (\ref{(bifix1)})).

We obtain a Sheffer sequence $\left( s_{n}\right) $, say, but only weakly
above the line $y=x+\delta $, 
\begin{equation*}
D\left( n,m\right) =s_{n}\left( m\right) \text{\ for all\ }m\geq \delta ,
\end{equation*}%
if the correction terms $D\left( n-a,m-c\right) $ (or $\sum\limits_{i\geq
0}(-1)^{i}D(n-a-bi,m-c-di)$) in (\ref{(bifix0')}) or (\ref{(bifix1')}) are
taken at points weakly above the same line, thus $a\geq c$ and $b\geq d$.
Finally, we must enforce in case of a nonempty bifix $o$ that $\func{depth}%
\left( op^{\prime }\right) $, where $p=op^{\prime }o$, does not exceed $%
d+a-c $.

Looking closer at the boundary $m=n+\delta $,%
\begin{equation*}
s_{n}\left( n+\delta \right) =s_{n-1}(n+\delta )+s_{n}(n+\delta
-1)-\sum\limits_{i\geq 0}(-1)^{i}s_{n-a-bi}(n+\delta -c-di)
\end{equation*}%
we see that $s_{n}\left( n+\delta -1\right) $ and $D\left( n,n+\delta
-1\right) $ also agree. We utilized this fact in the previous Sections to
determine the initial values $s_{n}\left( n-1\right) =\delta _{0,n}$ when
the depth was $0$.

\begin{example}
\label{ExDepth2} The following Table shows the number of ballot paths
avoiding the pattern $rrruuurrruu$, which has depth $\delta =2$, $a-c=1$,
one bifix $rrruu$, $op\prime =rrruuu$ of depth $3$, and $b-d=0$.
\end{example}

\tabcolsep=3pt%
$%
\begin{tabular}{l}
$%
\begin{tabular}{l||lllllllllll}
$12$ & \textbf{1} & \textbf{12} & \textbf{77} & \textbf{350} & \textbf{1260}
& \textbf{3808} & \textbf{9991} & \textbf{23219} & \textbf{48304} & \textbf{%
90046} & \textbf{148452} \\ 
$11$ & \textbf{1} & \textbf{11} & \textbf{65} & \textbf{273} & \textbf{910}
& \textbf{2548} & \textbf{6184} & \textbf{13235} & \textbf{25112} & \textbf{%
41816}$_{\bullet }$ & \emph{58567} \\ 
$10$ & \textbf{1} & \textbf{10} & \textbf{54} & \textbf{208} & \textbf{637}
& \textbf{1638} & \textbf{3637} & \textbf{7057} & \textbf{11\thinspace 897}
& \emph{16751}$_{\bullet }^{\Uparrow }$ & 16751 \\ 
$9$ & \textbf{1} & \textbf{9} & \textbf{44} & \textbf{154} & \textbf{429} & 
\textbf{1001} & \textbf{2000}$_{\bullet \Rightarrow }$ & \textbf{3425}$%
_{\bullet \Rightarrow }$ & \emph{4854}$_{\bullet \Rightarrow }$ & 
\multicolumn{1}{r}{4854$_{\bullet }^{\Uparrow }$} & 0 \\ 
$8$ & \textbf{1} & \textbf{8} & \textbf{35} & \textbf{110} & \textbf{275} & 
\textbf{572} & \textbf{1000}$_{\bullet }^{\Uparrow }$ & \emph{1429} & 1429 & 
0 &  \\ 
$7$ & \textbf{1} & \textbf{7} & \textbf{27} & \textbf{75} & \textbf{165} & 
\textbf{297} & \multicolumn{1}{r}{\emph{429}$_{\bullet }^{\uparrow \Uparrow
}~$} & 429 & 0 &  &  \\ 
$6$ & \textbf{1} & \textbf{6} & \textbf{20} & \textbf{48}$_{\bullet
\rightarrow }^{\hspace{4pt}\Rightarrow }$ & \textbf{90}$_{\bullet
\rightarrow }^{\hspace{4pt}\Rightarrow }$ & \emph{132}$_{\bullet \rightarrow
}^{\hspace{4pt}\Rightarrow }$ & \multicolumn{1}{r}{132$_{\bullet }^{\uparrow
\Uparrow }~$} & 0 &  &  &  \\ 
$5$ & \textbf{1} & \textbf{5} & \textbf{14} & \textbf{28}$_{\bullet
}^{\uparrow \boldsymbol{\Uparrow }}$ & \emph{42} & 42 & 0 &  & 
\multicolumn{3}{l}{The numbers in \emph{italics}} \\ 
$4$ & \textbf{1} & \textbf{4} & \textbf{9} & \emph{14}$_{\bullet }^{\uparrow 
\boldsymbol{\Uparrow }}$ & 14 & 0 &  &  & \multicolumn{3}{l}{are both,
values of the} \\ 
$3$ & \textbf{1}$_{\bullet \rightarrow }^{\hspace{4pt}\boldsymbol{%
\Rightarrow }}$ & \textbf{3}$_{\bullet \rightarrow }^{\hspace{4pt}%
\boldsymbol{\Rightarrow }}$ & \emph{5}$_{\bullet \rightarrow }^{\hspace{4pt}%
\boldsymbol{\Rightarrow }}$ & \multicolumn{1}{r}{5$_{\bullet }^{\uparrow 
\boldsymbol{\Uparrow }}\ $} & 0 &  &  &  & \multicolumn{3}{l}{Sheffer
sequence $\left( s_{n}\right) $} \\ 
$2$ & \textbf{1} & \emph{2} & 2 & 0 &  &  &  &  & \multicolumn{3}{l}{and of
the ballot recursion} \\ 
$1$ & \emph{1} & 1 & 0 &  &  &  &  &  &  &  &  \\ 
$0$ & 1 & 0 &  &  &  &  &  &  &  &  &  \\ \hline\hline
& $0$ & $1$ & $2$ & $3$ & $4$ & $5$ & $6$ & $7$ & $8$ & $9$ & $n$%
\end{tabular}%
$ \\ 
The path ($\rightarrow $) to $\left( 6,8\right) $ shows the first instance
when the occurrence of the pattern is \\ 
subtracted. The double path ($\Rightarrow $) to $\left( 9,11\right) $ shows
the first instance when the occurrence \\ 
of the $op^{\prime }$ path is added.%
\end{tabular}%
$

How can we determine initial values for the Sheffer sequence $\left(
s_{n}\right) $ if $\delta >0$? There are several ways of doing this; one way
applies the following Lemma to obtain a (recursive) formula for the values $%
s_{n}\left( n+\delta -1\right) $. Such a recursion will be enough for
explicitly representing the Sheffer polynomials in terms of their basic
sequence.

\begin{lemma}
\label{LemInitials}Let $a=\#r$'s in $p$ be larger than $1$. For all depths $%
\delta \geq 0$ holds $s_{0}\left( x\right) =1$, and for $0\leq j\leq \delta $
holds%
\begin{equation*}
D\left( n,n+j-1\right) =\sum_{i=1}^{j/2}\dbinom{j-i}{i}\left( -1\right)
^{i-1}D\left( n-i,n-i+j-1\right) +\sum_{i=1}^{\left( j+1\right) /2}\dbinom{%
j-i}{i-1}\left( -1\right) ^{i-1}D\left( n-i,n-i+j\right)
\end{equation*}%
for all $n\geq 1$. Especially%
\begin{equation*}
s_{n}\left( n+\delta -1\right) =\sum_{i=1}^{\delta /2}\dbinom{\delta -i}{i}%
\left( -1\right) ^{i-1}s_{n-i}\left( n-i+\delta -1\right)
+\sum_{i=1}^{\left( \delta +1\right) /2}\dbinom{\delta -i}{i-1}\left(
-1\right) ^{i-1}s_{n-i}\left( n-i+\delta \right)
\end{equation*}
\end{lemma}

\begin{proof}
For $n=1$ we get%
\begin{equation*}
D\left( 1,j\right) =\dbinom{j-1}{1}D\left( 0,j-1\right) +\dbinom{j-1}{0}%
D\left( 0,j\right) =j
\end{equation*}%
which is right because $a>1$. For $j=0$ we get $D\left( n,n-1\right) =0$ as
required. The double induction over $n\geq 1$ and $1\leq j\leq \delta $ uses
only the recursion $D\left( n,n-1+j\right) =D\left( n-1,n-1+j\right)
+D\left( n,n-2+j\right) $.
\end{proof}

\begin{example}[continued]
In the example avoiding the pattern $rrruuurrruu$, which has $\delta =2$, we
find $s_{n}\left( n+1\right) =s_{n-1}\left( n\right) +s_{n-1}\left(
n+1\right) $.
\end{example}

Using the above initial values from Lemma \ref{LemInitials} to find $%
s_{n}\left( x\right) $ above the line $y=x+\delta $ requires the Functional
Expansion Theorem (\cite[Theorem 3]{Nied03}):

\begin{theorem}
Suppose $\left( s_{n}\right) $ is a Sheffer sequence for the delta operator $%
B$ with basic sequence $\left( b_{n}\right) $, and $L$ a functional such
that $\left\langle L\mid 1\right\rangle \neq 0$. Then%
\begin{equation*}
s_{n}\left( x\right) =\sum_{k=0}^{n}\left\langle L\mid s_{k}\right\rangle 
\tilde{L}^{-1}b_{n-k}\left( x\right)
\end{equation*}%
where $\tilde{L}^{-1}$ is the operator inverse to $\tilde{L}=\sum_{i\geq
0}\left\langle L\mid b_{i}\right\rangle B^{i}$.
\end{theorem}

Lemma \ref{LemInitials} tells us how to choose the functional $L$,%
\begin{eqnarray*}
&&\left\langle L\mid s_{n}\left( n+x\right) \right\rangle \\
&=&s_{n}\left( n+\delta -1\right) -\sum_{i=1}^{\delta /2}\dbinom{\delta -i}{i%
}\left( -1\right) ^{i-1}s_{n-i}\left( n-i+\delta -1\right)
-\sum_{i=1}^{\left( \delta +1\right) /2}\dbinom{\delta -i}{i-1}\left(
-1\right) ^{i-1}s_{n-i}\left( n-i+\delta \right) \\
&=&\delta _{0,n}
\end{eqnarray*}%
Note that we applied the functional $L$ to $s_{n}\left( n+x\right) $, a
Sheffer polynomial for the delta operator $B_{1}:=E^{-1}B$, with basic
polynomials $\left( x+n\right) b_{n}\left( x\right) /x$ (see \cite[Remark 1]%
{Nied03}). In terms of the evaluation functional $\limfunc{Eval}%
\nolimits_{z}f\left( x\right) =f\left( z\right) $ we can write%
\begin{eqnarray*}
\left\langle L\mid s_{n}\left( n+x\right) \right\rangle &=&\left( \limfunc{%
Eval}\nolimits_{\delta -1}-\limfunc{Eval}\nolimits_{\delta
-1}\sum_{i=1}^{\delta /2}\dbinom{\delta -i}{i}\left( -1\right)
^{i-1}B_{1}^{i}\right. \\
&&\left. -\limfunc{Eval}\nolimits_{\delta }\sum_{i=1}^{\left( \delta
+1\right) /2}\dbinom{\delta -i}{i-1}\left( -1\right) ^{i-1}B_{1}^{i}\right)
s_{n}\left( n+x\right)
\end{eqnarray*}%
According to \cite[(2.7)]{Nied03} we get%
\begin{equation*}
\tilde{L}=E^{\delta -1}-E^{\delta -1}\sum_{i=1}^{\delta /2}\dbinom{\delta -i%
}{i}\left( -1\right) ^{i-1}B_{1}^{i}-E^{\delta }\sum_{i=1}^{\left( \delta
+1\right) /2}\dbinom{\delta -i}{i-1}\left( -1\right) ^{i-1}B_{1}^{i}
\end{equation*}%
and finally%
\begin{equation*}
\tilde{L}^{-1}=E^{1-\delta }/\left( 1-\sum_{i=1}^{\left( \delta +1\right)
/2}\left( \dbinom{\delta -i}{i}+E^{1}\dbinom{\delta -i}{i-1}\right) \left(
-1\right) ^{i-1}B_{1}^{i}\right)
\end{equation*}%
Therefore, the Functional Expansion Theorem writes $s_{n}\left( x\right) $
in terms of $b_{n}\left( x\right) $ as

$s_{n}\left( n+x\right) =$%
\begin{eqnarray*}
&&\tilde{L}^{-1}\frac{x}{n+x}b_{n}\left( n+x\right)  \\
&=&\frac{E^{1-\delta }}{1+\sum_{i=1}^{\left( \delta +1\right) /2}\left( 
\dbinom{\delta -i}{i}+E^{1}\dbinom{\delta -i}{i-1}\right) \left( -1\right)
^{i}B_{1}^{i}}\frac{x}{n+x}b_{n}\left( n+x\right)  \\
&=&\sum_{j\geq 0}\left( -1\right) ^{j}\left( \sum_{i=1}^{\left( \delta
+1\right) /2}\left( \dbinom{\delta -i}{i}+E^{1}\dbinom{\delta -i}{i-1}%
\right) \left( -1\right) ^{i}B_{1}^{i}\right) ^{j}\frac{x+1-\delta }{%
n+x+1-\delta }b_{n}\left( n+x+1-\delta \right) 
\end{eqnarray*}%
We leave the final expansion to the reader. Note that \newline
$B_{1}^{k}\frac{x+\delta -1}{n+x+\delta -1}b_{n}\left( n+x+1-\delta \right) =%
\frac{x+\delta -1}{n-k+x+\delta -1}b_{n-k}\left( n-k+x+1-\delta \right) $.
Of course, thevalues $s_{n}\left( n+x\right) $ agree with the pattern
avoiding path counts only for $x\geq \delta -1$.

\begin{example}[continued]
In the example above, avoiding the pattern $rrruuurrruu$, which has $\delta
=2$, we find 
\begin{eqnarray*}
s_{n}\left( n+x\right) &=&\sum_{j\geq 0}\left( 1+E^{1}\right) ^{j}B_{1}^{j}%
\frac{x-1}{n+x-1}b_{n}\left( n+x-1\right) \\
&=&2^{n}+\sum_{j=1}^{n}\sum_{i=0}^{n-j}\binom{n-j}{i}\frac{x+i-1}{j+x+i-1}%
b_{j}\left( j+x+i-1\right) .
\end{eqnarray*}
\end{example}

Finally, we also have to expand the ballot like numbers between the lines $%
y=x$ and $y=x+\delta $ in terms of those weakly above the line $y=x+\delta
-1 $.

\begin{lemma}
We have for all $0\leq j\leq \delta -1$
\end{lemma}

$D\left( n,n+\delta -1-j\right) =\sum_{i=0}^{j/2}\binom{j-i}{i}\left(
-1\right) ^{i}s_{n-i}\left( n-i+\delta -1\right) +\sum_{i=1}^{\left(
j+1\right) /2}\binom{j-i}{i-1}\left( -1\right) ^{i}s_{n-i}\left( n-i+\delta
\right) $.

\begin{proof}
Similar to the proof of Lemma \ref{LemInitials}.
\end{proof}

\begin{example}
In the above example, avoiding the pattern $rrruuurrruu$, we have $\delta =2$%
. We find $D\left( n,n+1-j\right) =s_{n}\left( n+1\right) -\binom{j-1}{0}%
s_{n-1}\left( n+1\right) $ for $j=1,2$. For example, $D\left( 10,10\right)
=s_{10}\left( 11\right) -s_{9}\left( 11\right) =58572-41821=\allowbreak
16751 $.
\end{example}

\section{The patterns of length 4}

In their paper \emph{Counting Strings in Dyck Paths }\cite{stt}, A.
Sapounakis, I. Tasoulas, and P. Tsikouras find generating functions for all
patterns of length 4 occurring $k$ times in Dyck paths, i.e., in ballot
paths \emph{returning to the diagonal}. Their case $k=0$ is our pattern
avoiding case. All patterns of length four are included in the above
considerations, except $uuur$, $uuru$, $uruu$, and $ruuu$. The path counts
in these four cases are not (eventually) Sheffer polynomial; the Finite
Operator Calculus does not directly apply.

\begin{description}
\item[$ruuu$] After differencing, this case is similar to $u^{3}$, and we
obtain for $m\geq 1$%
\begin{equation*}
D\left( n,n+m;ru^{3}\right) =\dbinom{n}{n}_{3}+\dbinom{n}{n-m}%
_{3}+2\sum_{i=1}^{m-1}\dbinom{n}{n-i}_{3}
\end{equation*}%
Note that $D\left( n,n;ru^{3}\right) =D\left( n-1,n;ru^{3}\right) $.

\item[$uuur$] In a similar way as in the case $ru^{3}$, we obtain%
\begin{equation*}
D\left( n,n+m;u^{3}r\right) =\sum_{i=0}^{m}\sum\limits_{k=0}^{i+1}\dfrac{1}{%
n+i+1-k}\dbinom{k-i-1}{k}_{3}\dbinom{n+i+1-k}{n+i-k}_{3}
\end{equation*}%
which also holds for $m=0$.

\item[$uuru$ and $uruu$] Both cases are avoided by the same ballot paths. We
obtain\newline
$D\left( n,n+l;uuru\right) =D\left( n,n+l;uruu\right) $%
\begin{eqnarray*}
&=&\sum_{k=0}^{n}\binom{m-1}{n-k}\left( -1\right) ^{n-k}\sum_{i=0}^{k\ /2}%
\dbinom{2i-k-1}{i}\frac{1}{k+1-2i}\times \\
&&\times \left( m\dbinom{2k+m-3i}{k-2i-1}+\dbinom{2k+m-3i}{k-2i}\right)
\end{eqnarray*}%
Note that except for the special cases $u^{4}$ and $r^{4}$ no pattern of
length $4$ has a bifix index larger than $1$.
\end{description}

\end{document}